\documentclass[a4paper]{article}

\usepackage[english]{babel}
\usepackage[utf8x]{inputenc}
\usepackage[T1]{fontenc}

\usepackage[top=3cm,bottom=2cm,left=3cm,right=3cm,marginparwidth=1.75cm]{geometry}
\usepackage{authblk}
\usepackage{blindtext}
\usepackage{amsmath, amsthm, amssymb, esint, mathtools, mathabx}
\usepackage{amsfonts}
\usepackage{graphicx}
\graphicspath{{/texfiles/}}
\usepackage[colorinlistoftodos]{todonotes}
\usepackage[colorlinks=true, allcolors=blue]{hyperref}
\usepackage{multicol}

\numberwithin{equation}{section}

\title{The pyramid averaging operator}

\author{A. Martina Neuman}
\affil{Department of Mathematics, New York University, Shanghai}
\affil[ ]{\textit{marsneuman@nyu.edu}}

\begin{document}

\maketitle

\begin{abstract} 

\noindent
This paper gives a concept of an integral operator defined on a manifold $M$ consisting of triple of points in $\mathbb{R}^{d}$ making up a regular $3$-simplex with the origin. The boundedness of such operator is investigated. The boundedness region contains more than the Banach range - a fact that mirrors the spherical $L^{p}$-improving estimate. The purpose of this paper is two-fold: one is to investigate into an integral operator over a manifold created from high-dimensional regular simplices, two is to start a maximal operator theory for such integral operator.

\end{abstract}

\tableofcontents

\section{Introduction}

\noindent
The so-called pyramid operator (pyramid averaging operator) belongs to a class of multilinear convolution surface averaging operators, ie, operators of the form:
\begin{equation}\label{eq I1}
(f_1,\cdots,f_{n})\mapsto [(f_1\otimes\cdots\otimes f_{n})\star\mu](x,\cdots,x)\end{equation}
where $f_{i}$'s are measurable functions on $\mathbb{R}^{d}$ and $\mu$ are a Borel measure on $\mathbb{R}^{nd}$. Examples of this type include the spherical averaging operator, the triangle averaging operator [7]. Operators of this type are typically used to investigate various continuous forms of discrete distance-set problems. For instance, the authors of [5] introduced a bilinear convolution operator 
\begin{equation*} B(f,g)(x)=\int_{M} f(x − u)g(x − v)\,dK(u, v)\end{equation*}
where $dK$ is the surface measure on the manifold $M=\{(u, v)\in\mathbb{R}^2\times\mathbb{R}^2: |u| = |v| = |u − v| = 1\}$. In other words, $M$ is the manifold consists of points $(u,v)$ such that $\{0, u, v\}$ forms an equilateral triangle in the plane. Then, the authors of [5] were able to show that, using their estimates on $B(f,g)$, if $E\subset\mathbb{R}^2$ is a compact set with $dim_{H}(E)>7/4$, then the set of three-point configurations determined by $E$ has a positive Lebesgue measure as a subset of $\mathbb{R}^3$.\\

\noindent
This paper is a generalization of the triangle operator dealt with in [7]. The author intends to initiate the study of the four-point configuration problem through a trilinear convolution operator - the pyramid averaging operator, as well as the study of the corresponding maximal operator, which can be thought of as a higher-dimensional version of the spherical maximal operator.\\

\noindent
In what follows, a "pyramid" means a $3$-simplex regular tetrahedron. Temporarily let $d\geq 4$. Let
\begin{equation*}M=\{(u,v,w)\in\mathbb{S}^{d-1}\times\mathbb{S}^{d-1}\times\mathbb{S}^{d-1}: |u-v|=|v-w|=|w-u|=1\}.\end{equation*}
In other words, $M$ consists of all points $(u,v,w)\in\mathbb{R}^{d}\times\mathbb{R}^{d}\times\mathbb{R}^{d}$ such that $\{0,u,v,w\}$ forms a pyramid. Let $\mu$ denote the surface measure on the manifold $M$. Define a trilinear pyramid operator $T$ as follows,
\begin{equation*}T(f,g,h)(x)=\int_{M}f(x-u)g(x-v)h(x-w)\,d\mu(u,v,w),\end{equation*}
with $f,g,h\in\mathcal{S}(\mathbb{R}^{d})$. One simply can't do anything with such generalization unless one transforms this definition into a more recognizable integration. To this end, let $(u,v,w)\in\mathbb{S}^{d-1}\times\mathbb{S}^{d-1}\times\mathbb{S}^{d-1}$, and
\begin{equation*}|u-v|^2=2-2\cos\theta(u,v);\,\,\,
|v-w|^2=2-2\cos\theta(v,w);\,\,\,|u-w|^2=2-2\cos\theta(u,w),\end{equation*}
where $\theta(x_1,x_2)$ denotes the angle between the two vectors $x_1,x_2$. Hence, in order for $(u,v,w)\in M$,
\begin{equation*}|\theta(u,v)|=|\theta(v,w)|=|\theta(u,w)|=\pi/3.\end{equation*}

\noindent
Consider one such $(u,v,w)\in M$. Choose $R\in SO(d)$ such that $u=Re_1, v=R(a_{v}e_1+b_{v}e_2), w=R(a_{w}e_1+b_{w}e_2+c_{w}e_3)$, with $a_{v}, b_{v}, a_{w}, b_{w}, c_{w}$ all being positive real numbers bounded by $1$ - they are sines and cosines of angles that make up the pyramid $\{0,u,v,w\}$. They are easily calculated and will be given in later sections. This description means that an integration over $M$ can be rewritten as an integration over $SO(d)$,
\begin{equation*}
\int_{M} F(x,y,z)\, d\mu(x,y,z) = \int_{SO(d)}F(Re_1,R(a_{v}e_1+b_{v}e_2),R(a_{w}e_1+b_{w}e_2+c_{w}e_3))\,dR.\end{equation*}

\noindent
Let $F(x,y,z) = f(x)g(y)h(z)$ with $f,g,h\in\mathcal{S}(\mathbb{R}^{d})$. Then a form of the pyramid operator is,
\begin{multline*}
T(f,g,h)(x) =\int_{M}f(x-u)g(x-v)h(x-w)\,d\mu(u,v,w)\\=\int_{SO(d)} f(x-Re_1)g(x- [Ra_{v}e_1+Rb_{v}e_2])\\
\times h(x-[Ra_{w}e_1+Rb_{w}e_2+Rc_{w}e_3])\,dR.\end{multline*}
That means that the Fourier transform of $\mu$ can be written as:
\begin{equation*}
\hat{\mu}(\xi,\delta,\eta)=\int_{SO(d)}\exp(-2\pi i\big[\xi\cdot Re_1+\delta\cdot R(a_{v}e_1+b_{v}e_2)\\
+\eta\cdot R(a_{w}e_1+(b_{w}e_2+c_{w}e_3\big])\,dR.\end{equation*}
The author investigates the boundedness of $T$ via the decay of $\hat{\mu}$. A tool needed for such analysis is {\bf Lemma A} given in {\bf Section 4}. However, in order to make use of the lemma, one needs to further decomposes $m$ into pieces $m_{i}$ that are better integrable with better controlled supports. See $\ref{eq E2}$. The price for such decomposition is the need for increase in dimensions. The result here is that, if $d>15$, then the boundedness region for $T$ contains more than just the Banach range. See {\bf Section 5} for a brief discussion of this range. This is the same phenomenon encountered in [7], where the triangle operator was considered. In [7], the requirement is $d>3$. This is a feature of this style of analysis. If one wants to lift this dimensional restriction, one might want to consider a different analysis route.

\subsection{Notations explanation}

\noindent
As usual, $\mathcal{S}(\mathbb{R}^{d})$ denotes the set of Schwartz functions on $\mathbb{R}^{d}$; $Conv(S)$ denotes the interior of the convex hull formed by points in $S$; $|\cdot|$ denotes either an absolute, a vector norm or a full-dimensional Lebesgue measure of a set; $\mathbb{S}^{d-1}$ denotes the unit sphere in $\mathbb{R}^{d}$ and $\mathbb{B}_{d}$ the unit ball in $\mathbb{R}^{d}$.

\section{Main theorem}

\noindent
{\bf Theorem 1.} Let the pyramid operator $T(f,g,h)(x)$ defined as in {\bf Section 1}. If $d>15$ then $T: L^{p}\times L^{q}\times L^{s}\to L^{r}$ in the region $(1/p,1/q,1/s)\in Conv(S)$ where
\begin{equation*}S=\{(0,1,0),(1,0,0),(0,0,1),(1/2,1/2,1/2),(1/p_0,1/p_0,0),(1/p_0,0,1/p_0),(0,1/p_0,1/p_0)\}\end{equation*}
with $p_0=\frac{5d}{3d-2}$ and $1/r =1/p+1/q+1/s$.

\section{Analytic tools needed}

\subsection{The slicing formula for spherical integral}

\begin{align*}
\nonumber \int_{\mathbb{S}^{d-1}} F(u)\, du  &=\sum_{\pm}\int_{\mathbb{B}^{d-1}} F(\pm\sqrt{1-|y|^2},y)\,\frac{dy}{\sqrt{1-|y|^2}} \\ &=\sum_{\pm}\int_0^1\int_{\mathbb{S}^{d-2}}F(\pm\sqrt{1-r^2},r\omega)\frac{r^{d-2}}{\sqrt{1-r^2}}\,d\omega dr.\end{align*}

\subsection{Quotient integral formula} 

\noindent
Suppose one has the following integration: $\int_{SO(d)}f(R)\,dR$, where $SO(d)$ is the orthogonal group in dimension $d$ and $dR$ is the left- and right-invariant Haar measure on $SO(d)$. For any closed subgroup $\mathcal{H}$ of $SO(d)$, there is an invariant Radon measure $d[S]$ on the quotient $SO(d)/\mathcal{H}$ such that,
\begin{equation*}\int_{SO(d)}f(R)\,dR = \int_{SO(d)/\mathcal{H}}\int_{\mathcal{H}}f(SR')\,dR'd[S].\end{equation*}
This is the content of {\bf Theorem 1.53} in [2]. In particular, one can normalize the measures to have them as probability measures. The map $S\mapsto\int_{\mathcal{H}} f(SR')\,dR'$ is then constant on cosets in $SO(d)/\mathcal{H}$, and one further has,
\begin{equation*}\int_{SO(d)/\mathcal{H}}\int_{\mathcal{H}} f(SR')\,dR'd[S]=\int_{SO(d)/\mathcal{H}}\int_{\mathcal{H}}\int_{\mathcal{H}} f(STR')\,dR'dTd[S] = \int_{SO(d)}\int_{\mathcal{H}} f(RR')\,dR'dR.\end{equation*}
For example, if $\mathcal{H}$ is the closed subgroup of rotations of fixing one vector $e_1$, then each coset $R\mathcal{H}$ consists of all rotations that maps $e_1\mapsto Re_1$ and $\mathcal{H}$ is isomorphic to $SO(d-1)$. One has the following final form of the quotient integral formula that is frequently used in the paper:
\begin{equation*}\int_{SO(d)}f(R)\,dR = \int_{SO(d)}\int_{SO(d-1)} f(RR')\,dR'dR.\end{equation*}

\subsection{Bessel functions}

\noindent
Let $J_{s}$ denote the Bessel function of order $s$. For a reference text on the following facts for $Re(s)>-1/2$, see [9].

\subsubsection{Recurrence formula}

\noindent
\begin{equation*}
\frac{d}{dt}(t^{-s}J_{s}(t))=-t^{-s}J_{s+1}(t).\end{equation*}

\subsubsection{Bessel functions of small and large arguments}

\noindent
Suppose $t\to 0^{+}$. Then,
\begin{equation*}J_{s}(t) = \frac{t^{s}}{2^{s}\Gamma(s+1)}+S_{s}(t)\end{equation*}
with 
\begin{equation*}|S_{s}(t)|\leq\frac{2^{-Re(s)}t^{Re(s)+1}}{(Re(s)+1)|\Gamma(s+1/2)|\Gamma(1/2)}.\end{equation*}

\noindent
In particular this means that if $s>0$ then as $t\to 0^{+}$, $\big|\frac{J_{s}(t)}{t^{s}}\big|\lesssim_{s} 1$.\\

\noindent
For $t\geq 1$, one has, $|J_{s}(t)|\lesssim_{s} t^{-1/2}$; the dominant constant can be chosen to depend smoothly on $Re(s)\in (-1/2,\infty)$.\\

\noindent
{\it Remark 1:} These facts, the recurrence formula and the growth rates, imply that if $\alpha>0$, $|(d^{\alpha}/dt^{\alpha})(t^{-s}J_{s}t)|\lesssim_{\alpha,s} 1$ as $t\to 0^{+}$, and, $|(d^{\alpha}/dt^{\alpha})(t^{-s}J_{s}(t))|\approx_{\alpha,s} o(t^{-s})+t^{-s}J_{u}(t)$ for some $u=u(\alpha)>0$, as $t\to\infty$. Altogether, one has,
\begin{equation}\label{eq A1}
\bigg|\frac{d^{\alpha}}{dt^{\alpha}}(t^{-s}J_{s}(t))\bigg|\leq O((1+t)^{-s})|J_{u}(t)|.\end{equation}
Here $O((1+t)^{-s})$ denotes a term whose magnitude is of the specified order. 

\subsubsection{The Fourier transform of surface measure}

\noindent
Let $d\sigma$ be the surface measure on $\mathbb{S}^{d-1}$, $d\geq 2$. Then,
\begin{equation}\label{eq A2}
\hat{\sigma}_{d-1}(\xi)=\int_{\mathbb{S}^{d-1}} \exp(-2\pi i\xi\cdot\theta)\,d\theta = \frac{2\pi J_{\frac{d-2}{2}}(2\pi|\xi|)}{|\xi|^{\frac{d-2}{2}}}.\end{equation}

\section{Lemma needed}

\noindent
The following lemma, proven in [6], is needed for the subsequent analysis of the decomposition of $\hat{\mu}=:m$.\\

\noindent
{\bf Lemma A.} Let $1\leq q<3$ and set $M_{q}=\lfloor\frac{6d}{3-q}\rfloor +1$. Let $m(\xi,\delta,\eta)$ be a function in $L^{q}(\mathbb{R}^{3d})\cap\mathcal{C}^{M_{q}}(\mathbb{R}^{3d})$ satisfying
\begin{equation*} \|\partial^{\alpha}m\|_{L^{\infty}}\leq C_0 \end{equation*}
for all $|\alpha|\leq M_{q}$. Then the trilinear operator $T_{m}$ with the multiliplier $m$ satisfies 
\begin{equation*}
\|T_{m}\|_{L^2\times L^2\times L^2\to L^{2/3}}\lesssim_{C_0,d,q}\|m\|^{q/3}_{L^{q}}.\end{equation*}

\noindent
{\it Remark 2:} This statement of {\bf Lemma A} doesn't specify clearly the dependence of $\|T_{m}\|_{L^2\times L^2\times L^2\to L^{2/3}}$ on $C_0$. However a remark followed from the proof of {\bf Lemma A} in [6] made clear that one can take the dominant constant to be $C_0^{1-q/3}$. Now, as $1\leq q<3$, $d\geq 4$, one can take $3d =: M_{q}\geq\lfloor\frac{3d}{2}\rfloor+1$. Finally, when $m$ is a measure of finite support, then $\|m\|^{q/3}_{L^{q}}=|supp(m)|^{1/3}$.

\section{Banach range}

\noindent
The Banach range [4] for the boundedness of $T: L^{p}\times L^{q}\times L^{s}\to L^{r}$ is the convex hull 
\begin{equation*} Conv((1,0,0),(0,1,0),(0,0,1))\ni (1/p,1/q,1/s) \end{equation*} 
with $1/r=1/p+1/q+1/s\leq 1$. This Banach range can be easily obtained by through an associated spherical integration. To see this, first recall that,
\begin{multline*}
T(f,g,h)(x) =\int_{M}f(x-u)g(x-v)h(x-w)\,d\mu(u,v,w)\\=\int_{SO(d)} f(x-Re_1)g(x- [Ra_{v}e_1+Rb_{v}e_2])\\
\times h(x-[Ra_{w}e_1+Rb_{w}e_2+Rc_{w}e_3])\,dR.\end{multline*}
That means through a crude domination, one arrives at,
\begin{equation*}
|T(f,g,h)(x)|\lesssim_{d}S_1(|f|)(x)\|g\|_{\infty}\|h\|_{\infty},\end{equation*}
where $S_1$ denotes the spherical average operator whose range of boundedness $S_1: L^{p}\to L^{r}$ is established within the region that is the interior of the following convex hull:
\begin{equation*} (1/p,1/r)\in V = Conv((0,0),(1,1),(d/(d+1),1/(d+1)).\end{equation*}
That means that $T: L^{p}\times L^{\infty}\times L^{\infty}\to L^{r}$ with $(1/p,1/r)\in V$. Due to symmetry between $f,g,h$ one then has:
\begin{align*} T &: L^{p}\times L^{\infty}\times L^{\infty}\to L^{r}\\
T &: L^{\infty}\times L^{\infty}\times L^{p}\to L^{r}\\
T &: L^{\infty}\times L^{p}\times L^{\infty}\to L^{r}
\end{align*}
with $(1/p,1/r)\in V$. Multilinear interpolation allows one to obtain the Banach range for $T$: 
\begin{equation*}T: L^{p}\times L^{q}\times L^{s}\to L^{r}\end{equation*}
with $(1/p,1/q,1/s)\in Conv((1,0,0),(0,1,0),(0,0,1))$ and $1/r = 1/p+1/q+1/s$, plus a little bit more due to the $L^{p}$-improving property of the spherical average operator. However, one can consider a different route to obtain a bigger range. For that, see {\bf Section 10}.

\section{A spherical calculation}

\noindent
Pick $(u,v,w)\in M$. Then one can let:
\begin{equation}\label{eq S1} u=Re_1,v=R((1/2)e_1+(\sqrt{3}/2)e_2), w=R((1/2)e_1+(1/2\sqrt{3})e_2+(\sqrt{2}/\sqrt{3})e_3),\end{equation} 
and an integration over $M$ is the same as an integration over $SO(d)$:
\begin{multline*}\int_{M}f(u,v,w)\,d\mu(u,v,w) = \int_{SO(d)} f(Re_1,R((1/2)e_1+R((1/2)e_1+(\sqrt{3}/2)e_2),\\
R((1/2)e_1+(1/2\sqrt{3})e_2+(\sqrt{2}/\sqrt{3})e_3))\,dR.\end{multline*}
Consider the integration that characterizes the Fourier transform of $\mu$:
\begin{multline}\label{eq S2}
\hat{\mu}(\xi,\delta,\eta)=\int_{SO(d)}\exp(-2\pi i\big[\xi\cdot Re_1+\delta\cdot R((1/2)e_1+(\sqrt{3}/2)e_2)\\
+\eta\cdot R((1/2)e_1+(1/2\sqrt{3})e_2+(\sqrt{2}/\sqrt{3})e_3\big])\,dR\end{multline}

\subsection{First decomposition}

\noindent
One utilizes the quotient integral first time here, with $\mathcal{H}\cong SO(d-1)$ being a subgroup of $SO(d)$ that fixes $e_1$. Then $\ref{eq S2}$ becomes:
\begin{multline}\label{eq S3} \hat{\mu}(\xi,\delta,\eta)=\int_{SO(d)}\exp(-2\pi i\big[\xi\cdot Re_1+(1/2)\delta\cdot Re_1+(1/2)\eta\cdot Re_1\big]
\\ \times\big(\int_{SO(d-1)}\exp(-2\pi i\big[(\sqrt{3}/2)\delta\cdot RR'e_2+(1/2\sqrt{3})\eta\cdot RR'e_2+(\sqrt{2}/\sqrt{3})\eta\cdot RR'e_3\big])\,dR'\big) dR.\end{multline}

\noindent
The inner integral of $\ref{eq S3}$ is:
\begin{equation*}\int_{SO(d-1)}\exp(-2\pi i\big[(\sqrt{3}/2)\delta\cdot RR'e_2+(1/2\sqrt{3})\eta\cdot RR'e_2+(\sqrt{2}/\sqrt{3})\eta\cdot RR'e_3\big] )\,dR'\end{equation*}
Utilizing the invariance of the Haar measure on $SO(d)$, one can precompose $R\in SO(d)$ with a rotation whose transpose maps
\begin{equation}\label{eq S4} \eta\mapsto|\eta|e_2,\delta\mapsto|\delta|a_2e_2+|\delta|a_3e_3,\xi\mapsto|\xi|b_1e_1+|\xi|b_2e_2+|\xi|b_3e_3,\end{equation} 
with $a_2=\cos\theta(\delta,e_2)=\cos\theta(\delta,\eta)$, $a_3=\cos\theta(\delta,e_3)=\sin\theta(\delta,\eta)$, $b_1=\cos\theta(\xi,e_1)$, $b_2=\cos\theta(\xi,e_2)=\cos\theta(\xi,\eta)$, $b_3=\cos\theta(\xi,e_3)$. This precomposition has an effect of rotating the sphere under consideration to a standard position where $\eta,\delta$ are both on the $yz$-plane. It also has an effect in transforming the inner integral of $\ref{eq S3}$:
\begin{equation}\label{eq S5} \int_{SO(d-1)}\exp(-2\pi i\big[(\sqrt{3}/2)(|\delta|a_2e_2+|\delta|a_3e_3)\cdot RR'e_2+(1/2\sqrt{3})|\eta|e_2\cdot RR'e_2+(\sqrt{2}/\sqrt{3})|\eta|e_2\cdot RR'e_3\big])\,dR'.\end{equation}
Using orthogonality and the fact that $R'e_2,R'e_3$ are two orthonormal unit vectors in $\mathbb{S}^{d-2}$, $\ref{eq S5}$ can be transformed into
\begin{multline}\label{eq S6}
\int_{SO(d-1)}\exp(-2\pi i\big[(\sqrt{3}/2) P_{-1}R^{T}(|\delta|a_2e_2+|\delta|a_3e_3)\cdot R'e_2+(1/2\sqrt{3})P_{-1}R^{T}|\eta|e_2\cdot R'e_2
\\ +(\sqrt{2}/\sqrt{3})P_{-1}R^{T}|\eta|e_2\cdot R'e_3\big])\,dR'.\end{multline}
Here, $P_{-1}$ is the projection map that maps every vector into its last $d-1$ coordinates, for example:
\begin{equation*}P_{-1}(1,2,3)=(2,3).\end{equation*}
The meaning of $P_{-1}$ changes with each dimension under consideration, but it always maps a vector in $\mathbb{R}^{d}$ to a vector in $\mathbb{R}^{d-1}$. 

\subsection{Second decomposition}

\noindent
Consider a subgroup $\mathcal{G}\cong SO(d-2)$ of $SO(d-1)$ that fixes $e_2$. Then use the quotient integral formula once more with $\mathcal{G}$ transform $\ref{eq S6}$ into:
\begin{multline}\label{eq S7}
\int_{SO(d-1)}\exp(-2\pi i\big[(\sqrt{3}/2) P_{-1}R^{T}(|\delta|a_2e_2+|\delta|a_3e_3)\cdot R'e_2+(1/2\sqrt{3})P_{-1}R^{T}|\eta|e_2\cdot R'e_2])
\\ \times\big(\int_{SO(d-2)}\exp(-2\pi i[(\sqrt{2}/\sqrt{3})P_{-1}R^{T}|\eta|e_2\cdot R'R''e_3\big])\,dR''\big)dR'\end{multline}
whose inner integral is,
\begin{equation}\label{eq S8}
\int_{SO(d-2)}\exp(-2\pi i[(\sqrt{2}/\sqrt{3})P_{-1}R^{T}|\eta|e_2\cdot R'R''e_3\big])\,dR''=\hat{\sigma}_{d-3}(|\eta|(\sqrt{2}/\sqrt{3})P_{-1}(R')^{T}P_{-1}R^{T}e_2).\end{equation}

\noindent
Now $P_{-1}R^{T}e_2$ is an $\mathbb{R}^{d-1}$, as explained about the function of $P_{-1}$. Hence $P_{-1}(R')^{T}P_{-1}R^{T}e_2$ is an $\mathbb{R}^{d-2}$ vector. That means, through Pythagorean theorem:
\begin{equation}\label{eq S9}
|P_{-1}(R')^{T}P_{-1}R^{T}e_2|^2=|(R')^{T}P_{-1}R^{T}e_2|^2-((R')^{T}P_{-1}R^{T}e_2\cdot j_1)^2=|P_{-1}R^{T}e_2|^2-(P_{-1}R^{T}e_2\cdot R'j_1)^2.\end{equation}
Here $j_1$ denotes the first canonical unit vector of $\mathbb{R}^{d-1}$ - the one that has the first component to be $1$ and the remaining $d-2$ components zero. Now put $\ref{eq S7},\ref{eq S8},\ref{eq S9}$ into $\ref{eq S6}$ to have:
\begin{multline}\label{eq S10}
\int_{SO(d-1)}\exp(-2\pi i\big[(\sqrt{3}/2) P_{-1}R^{T}(|\delta|a_2e_2+|\delta|a_3e_3)\cdot R'j_1+(1/2\sqrt{3})P_{-1}R^{T}|\eta|e_2\cdot R'j_1\big])\\
\times\hat{\sigma}_{d-3}(|\eta|(\sqrt{2}/\sqrt{3})\sqrt{|P_{-1}R^{T}e_2|^2-(P_{-1}R^{T}e_2\cdot R'j_1)^2}\big)\,dR'\end{multline}

\subsection{Analysis of $\ref{eq S10}$}

\noindent
Firstly, if $R^{T}$ is a $d\times d$ rotation matrix, then $R^{T}e_2, R^{T}e_3$ are respectively the second and third columns of $R^{T}$. These columns are orthogonal to each other. That means that 
\begin{multline*}0=(R^{T}e_2\cdot e_1)(R^{T}e_3\cdot e_1)+(R^{T}e_2\cdot e_2)(R^{T}e_3\cdot e_2)+\cdots+(R^{T}e_2\cdot e_{d})(R^{T}e_3\cdot e_{d})\\
\Rightarrow (R^{T}e_2\cdot e_2)(R^{T}e_3\cdot e_2)+\cdots+(R^{T}e_2\cdot e_{d})(R^{T}e_3\cdot e_{d})=-(R^{T}e_2\cdot e_1)(R^{T}e_3\cdot e_1).\end{multline*}
This means that if $M_{R},N_{R}$ denotes the norms of $P_{-1}R^{T}e_2,P_{-1}R^{T}e_3$, respectively, and $\theta'$ is the angle between $P_{-1}R^{T}e_2,P_{-1}R^{T}e_3$ then $\cos\theta'$ can be captured by
\begin{equation}\label{eq S11}
\cos\theta'=\frac{-(R^{T}e_2\cdot e_1)(R^{T}e_3\cdot e_1)}{M_{R}N_{R}}=\frac{-(e_2\cdot Re_1)(e_3\cdot Re_1)}{M_{R}N_{R}}.\end{equation}
Furthermore, 
\begin{align}\label{eq S12}
\nonumber M_{R}^2 &=|P_{-1}R^{T}e_2|^2=1-(R^{T}e_2\cdot e_1)^2=1-(e_2\cdot Re_1)^2\\
N_{R}^2 &=|P_{-1}R^{T}e_3|^2=1-(R^{T}e_3\cdot e_1)^2=1-(e_3\cdot Re_1)^2.\end{align}

\noindent
Utilizing the invariance of the Haar measure on $SO(d-1)$ again, if one precomposes $R'\in SO(d-1)$ with a rotation in $SO(d-1)$ whose transpose maps $P_{-1}R^{T}e_2\mapsto M_{R}e_2$ and $P_{-1}R^{T}e_3\mapsto N_{R}\cos\theta'e_2+N_{R}\sin\theta'e_3$, then $\ref{eq S10}$ becomes:
\begin{multline}\label{eq S13}
\int_{SO(d-1)}\exp(-2\pi i[(\sqrt{3}/2) M_{R}|\delta|a_2 +(\sqrt{3}/2) N_{R}|\delta|a_3\cos\theta'+(1/2\sqrt{3})M_{R}|\eta|)(e_2\cdot R'j_1) \\
+ (\sqrt{3}/2)N_{R}|\delta|a_3\sin\theta' (e_3\cdot R'j_1)]) \times\hat{\sigma}_{d-3}(|\eta|(\sqrt{2}/\sqrt{3})\sqrt{M_{R}^2-(M_{R}e_2\cdot R'j_1)^2})\,dR'.\end{multline}
Now $R'j_1$ is a unit vector in $\mathbb{R}^{d-1}$, that means one can rewrite $\ref{eq S13}$ in terms of a spherical integral:
\begin{multline}\label{eq S14}
\int_{\mathbb{S}^{d-2}}\exp(-2\pi i\big[(\sqrt{3}/2) M_{R}|\delta|a_2+(\sqrt{3}/2) N_{R}|\delta|a_3\cos\theta'+(1/2\sqrt{3})M_{R}|\eta|)(e_2\cdot v) \\+ (\sqrt{3}/2)N_{R}|\delta|a_3\sin\theta' (e_3\cdot v)\big])\times\hat{\sigma}_{d-3}(|\eta|(\sqrt{2}/\sqrt{3})\sqrt{M_{R}^2-(M_{R}e_2\cdot v)^2})\,d\sigma(v)\end{multline}

\noindent
Recall here the slicing formula for the spherical integral adapted to $\mathbb{S}^{d-2}$:
\begin{multline*}\int_{\mathbb{S}^{d-2}}f(u)\,du = \sum_{\pm}\int_{B_{d-2}} f\left(\pm\sqrt{1-|y|^2},y\right)\,\frac{dy}{\sqrt{1-|y|^2}}\\
=\sum_{\pm}\int_0^1\int_{\mathbb{S}^{d-3}}f\left(\pm\sqrt{1-r^2},r\omega\right)\frac{r^{d-3}}{\sqrt{1-r^2}}\,d\omega dr.\end{multline*}
Applying this spherical integral to $\ref{eq S14}$ - utilizing the fact that $e_2,e_3$ are orthogonal to each other - gives it the following value:
\begin{multline*}
\sum_{\pm}\int_0^1\int_{\mathbb{S}^{d-3}}\exp(-2\pi i\big[(\sqrt{3}/2) M_{R}|\delta|a_2+(\sqrt{3}/2) N_{R}|\delta|a_3\cos\theta'+(1/2\sqrt{3})M_{R}|\eta|)\sqrt{1-r^2}\\
+(\sqrt{3}/2)N_{R}|\delta|a_3\sin\theta'(e_3\cdot r\omega)\big]) \times\hat{\sigma}_{d-3}(|\eta|(\sqrt{2}/\sqrt{3})M_{R}r)\frac{r^{d-3}}{\sqrt{1-r^2}}\,d\omega dr\end{multline*}
which, after an integration over $\mathbb{S}^{d-3}$, becomes:
\begin{multline}\label{eq S15}
2\int_0^1\cos(2\pi(\sqrt{3}/2) M_{R}|\delta|a_2+(\sqrt{3}/2) N_{R}|\delta|a_3\cos\theta'+(1/2\sqrt{3})M_{R}|\eta|)\sqrt{1-r^2})
\\ \times\hat{\sigma}_{d-3}((\sqrt{3}/2)N_{R}|\delta|a_3 r\sin\theta')\hat{\sigma}_{d-3}(|\eta|(\sqrt{2}/\sqrt{3})M_{R}r)\frac{r^{d-3}}{\sqrt{1-r^2}}\,d\omega dr\end{multline}

\subsection{Putting everything together}

\noindent
Recall that the outer integral of $\ref{eq S3}$ is
\begin{equation*}\int_{SO(d)}\exp(-2\pi i\big[(|\xi|b_1e_1+|\xi|b_2e_2+|\xi|b_3e_3)\cdot Re_1+(1/2)(|\delta|a_2e_2+|\delta|a_3e_3)\cdot Re_1+(1/2)|\eta|e_2\cdot Re_1])\,dR\end{equation*}
with the $a_{i},b_{i}$ values as in $\ref{eq S4}$. The inner integral of $\ref{eq S3}$ is precisely $\ref{eq S15}$. Use the invariance property of the Haar measure on $SO(d)$ again, this time, by precomposing $R\in SO(d)$ with a map whose transpose maps 
\begin{equation*}e_1\mapsto e_3, e_3\mapsto -e_1,\end{equation*}
and everything else remains the same. One might ask if this will interfere with the previous computations, but the precomposition is done before all of those computations. Think of this as follows. $SO(d)\ni R=S_1S_2S_3$, where $S_3\in SO(d)$, $S_1$ is the original decomposition, which after whose peeling $R\mapsto S_2S_3$, gives $\ref{eq S15},\ref{eq S16}$. Now one gives another peeling $S_2S_3\mapsto S_3$. The final value of the inner integral of $\ref{eq S3}$, which is $\ref{eq S15}$ only depends on $R$ through $M_{R},N_{R},\theta'$, whose values are easily changed with this precomposition. Moreover there is a symmetry among $u,v,w$ in $\ref{eq S1}$. With that, the outer integral of $\ref{eq S3}$ becomes, 
\begin{equation}\label{eq S16}\int_{SO(d)}\exp(-2\pi i\big[(|\xi|b_1e_3+|\xi|b_2e_2-|\xi|b_3e_1)\cdot Re_1+(1/2)(|\delta|a_2e_2-|\delta|a_3e_1)\cdot Re_1+(1/2)|\eta|e_2\cdot Re_1])\,dR.\end{equation}
The equations $\ref{eq S11},\ref{eq S12}$ then become:
\begin{equation}\label{eq S17}
\cos\theta'=\frac{(e_1\cdot Re_1)(e_2\cdot Re_1)}{M_{R}N_{R}},\end{equation}
and,
\begin{align}\label{eq S18}
\nonumber M_{R}^2 &=1-(e_2\cdot Re_1)^2\\
N_{R}^2 &=1-(e_1\cdot Re_1)^2.\end{align}

\noindent
Under this scheme, putting $\ref{eq S15},\ref{eq S16}$ together, one arrives at the following form of $\ref{eq S3}$:
\begin{multline}\label{eq S19}\int_{SO(d)}\exp(-2\pi i\big[(|\xi|b_1e_3+|\xi|b_2e_2-|\xi|b_3e_1)\cdot Re_1+(1/2)(|\delta|a_2e_2-|\delta|a_3e_1)\cdot Re_1+(1/2)|\eta|e_2\cdot Re_1])
\\ \times 2\int_0^1\cos(2\pi(\sqrt{3}/2) M_{R}|\delta|a_2+(\sqrt{3}/2) N_{R}|\delta|a_3\cos\theta'+(1/2\sqrt{3})M_{R}|\eta|)\sqrt{1-r^2})
\\ \times\hat{\sigma}_{d-3}((\sqrt{3}/2)N_{R}|\delta|a_3 r\sin\theta')\hat{\sigma}_{d-3}(|\eta|(\sqrt{2}/\sqrt{3})M_{R}r)\frac{r^{d-3}}{\sqrt{1-r^2}}\,dr\, dR.\end{multline}

\subsubsection{The slicing formula}

\noindent
If one applies the slicing formula for the spherical integral twice, one will obtain:
\begin{align}\label{eq S20}
\nonumber\int_{\mathbb{S}^{d-1}}f(u)\,du &= \sum_{\pm}\int_{B_{d-2}} f\left(\pm\sqrt{1-|y|^2},y\right)\,\frac{dy}{\sqrt{1-|y|^2}}\\
\nonumber &=\sum_{\pm}\int_0^1\int_{\mathbb{S}^{d-2}}f\left(\pm\sqrt{1-t^2},t\omega\right)\frac{r^{d-2}}{\sqrt{1-r^2}}\,d\omega dt\\
&=\sum_{\pm}\sum_{\pm}\int_0^1\int_0^1\int_{\mathbb{S}^{d-3}}f(\pm\sqrt{1-r^2},\pm t\sqrt{1-s^2},tsz)\,\frac{t^{d-2}}{\sqrt{1-t^2}}\frac{s^{d-3}}{\sqrt{1-s^2}}\,dz dtds.\end{align}

\noindent
Now note that in $\ref{eq S19}$, $Re_1$ is a vector in $\mathbb{S}^{d-1}$. That allows one to rewrite $\ref{eq S19}$ back to a spherical integration, remembering $\ref{eq S17},\ref{eq S18}$:
\begin{align}\label{eq S21}
\nonumber \hat{\mu}(\xi,\delta,\eta) =\int_{\mathbb{S}^{d-1}}\exp(&-2\pi i \big[(|\xi|b_1e_3+|\xi|b_2e_2-|\xi|b_3e_1)\cdot v+(1/2)(|\delta|a_2e_2-|\delta|a_3e_1)\cdot v+(1/2)|\eta|e_2\cdot v])\\
\nonumber &\times 2\int_0^1\cos\bigg(\bigg(2\pi(\sqrt{3}/2) \sqrt{1-(e_2\cdot v)^2}|\delta|a_2+(1/2\sqrt{3})\sqrt{1-(e_2\cdot v)^2}|\eta|\\
\nonumber &+(\sqrt{3}/2) \sqrt{1-(e_1\cdot v)^2}|\delta|a_3\frac{(e_2\cdot v)(e_1\cdot v)}{\sqrt{1-(e_1\cdot v)^2}\sqrt{1-(e_2\cdot v)^2}}\bigg)\sqrt{1-r^2}\bigg)\\
\nonumber &\times \hat{\sigma}_{d-3}\bigg((\sqrt{3}/2)\sqrt{1-(e_1\cdot v)^2}|\delta|a_3 r\sqrt{1-\frac{(e_2\cdot v)^2(e_1\cdot v)^2}{(1-(e_1\cdot v)^2)(1-(e_2\cdot v)^2)}}\bigg)\\
&\times \hat{\sigma}_{d-3}(|\eta|(\sqrt{2}/\sqrt{3})\sqrt{1-(e_2\cdot v)^2}r)\frac{r^{d-3}}{\sqrt{1-r^2}}\, dr\, d\sigma(v).\end{align}

\noindent
If one wish to apply $\ref{eq S20}$ to $\ref{eq S21}$ one needs to take orthogonality into consideration. In $\ref{eq S20}$, $z\in\mathbb{S}^{d-3}$ and $\sqrt{1-s^2}$ plays the role of the height of $\mathbb{S}^{d-2}$ above $B_{d-3}$. That means that $e_2\cdot z$ can be taken to be the height $\sqrt{1-s^2}$. If one moves onto the next level, and now $\omega\in\mathbb{S}^{d-2}$ then $\sqrt{1-t^2}$ similarly plays the role of the height of $\mathbb{S}^{d-1}$ over $B_{d-2}$. That means one can identify the following in $\ref{eq S21}$:
\begin{align}\label{eq S22}
\nonumber e_1\cdot v &=\sqrt{1-t^2}\\
\nonumber e_2\cdot v &=t\sqrt{1-s^2}\\
e_3\cdot v &= k_1\cdot tsz\end{align}
where $k_1$ here denotes the first canonical unit vector in $\mathbb{R}^{d-2}$.\\

\noindent
Applying $\ref{eq S20},\ref{eq S22}$ to $\ref{eq S21}$, one has:
\begin{align}\label{eq S23}
\nonumber \hat{\mu}(\xi,\delta,\eta) &= 2\sum_{\pm}\sum_{\pm}\int_0^1\int_0^1\int_0^1 \int_{\mathbb{S}^{d-3}}
\exp(-2\pi i \big[|\xi|b_1 ts(k_1\cdot z)+ |\xi|b_2t\sqrt{1-s^2} -|\xi|b_3\sqrt{1-t^2})\\
\nonumber &+(1/2)(|\delta|a_2t\sqrt{1-s^2}-|\delta|a_3\sqrt{1-t^2})+(1/2)|\eta|t\sqrt{1-s^2}])\\
\nonumber &\times \cos\bigg\{\bigg(2\pi(\sqrt{3}/2) \sqrt{1-t^2(1-s^2)}|\delta|a_2+(\sqrt{3}/2) t|\delta|a_3\frac{\sqrt{1-s^2}\sqrt{1-t^2}}{\sqrt{1-t^2(1-s^2)}} +(1/2\sqrt{3})\sqrt{1-t^2(1-s^2)}|\eta|\bigg)\\
\nonumber &\times\sqrt{1-r^2}\bigg\}\times \hat{\sigma}_{d-3}\bigg((\sqrt{3}/2) t|\delta|a_3 r\sqrt{1-\frac{(1-t^2)(1-s^2)}{1-t^2(1-s^2)}}\bigg)\times\hat{\sigma}_{d-3}(|\eta|(\sqrt{2}/\sqrt{3})\sqrt{1-t^2(1-s^2)}r)\\
&\frac{r^{d-3}}{\sqrt{1-r^2}}\frac{s^{d-3}}{\sqrt{1-s^2}}\frac{t^{d-2}}{\sqrt{1-t^2}}\,d\sigma(z)drdsdt.\end{align}

\noindent
One notes that in $\ref{eq S23}$ the integration over $\mathbb{S}^{d-3}$ can be evaluated first, which gives the following value for $\hat{\mu}(\xi,\delta,\eta)$
\begin{align}\label{eq S24}
\nonumber 8\int_0^1\int_0^1\int_0^1 &\cos\bigg(|\xi|b_2\sqrt{1-t^2}-|\xi|b_3\sqrt{1-t^2})+ (1/2)(|\delta|a_2t\sqrt{1-s^2}-|\delta|a_3\sqrt{1-t^2})+(1/2)|\eta|t\sqrt{1-s^2}\bigg)\\
\nonumber &\times\cos\bigg\{\bigg(2\pi(\sqrt{3}/2) \sqrt{1-t^2(1-s^2)}|\delta|a_2+(\sqrt{3}/2) t|\delta|a_3\frac{\sqrt{1-s^2}\sqrt{1-t^2}}{\sqrt{1-t^2(1-s^2)}} +(1/2\sqrt{3})\sqrt{1-t^2(1-s^2)}|\eta|\bigg)\\ 
\nonumber &\times\sqrt{1-r^2}\bigg\}\times\hat{\sigma}_{d-3}\bigg((\sqrt{3}/2) t|\delta|a_3 r\sqrt{1-\frac{(1-t^2)(1-s^2)}{1-t^2(1-s^2)}}\bigg)\times\hat{\sigma}_{d-3}(|\xi|b_1 ts)\\
&\times\hat{\sigma}_{d-3}(|\eta|(\sqrt{2}/\sqrt{3})\sqrt{1-t^2(1-s^2)}r)\frac{r^{d-3}}{\sqrt{1-r^2}}\frac{t^{d-2}}{\sqrt{1-t^2}}\frac{s^{d-3}}{\sqrt{1-s^2}}\,drdsdt.\end{align}
The only value among $a_{i},b_{i}$ in $\ref{eq S24}$ that one should remember, for the subsequent analysis is:
\begin{equation*} 
a_3 =\sin\theta(\delta,\eta),b_1=\cos\theta(\xi,e_3).\end{equation*}

\section{An estimation on derivatives}

\noindent
Expression $\ref{eq S24}$ might look intimidating, but one only needs to care about the non-cosine factors, as when one takes derivatives of $\ref{eq S24}$ in terms of $\xi,\delta,\eta$, the contribution of derivatives from those cosine factors are only at most a constant in absolute values. Hence the factors that one needs to consider is:
\begin{align}\label{eq E1}
\nonumber 8\int_0^1\int_0^1\int_0^1 &\hat{\sigma}_{d-3}\bigg((\sqrt{3}/2) t|\delta|a_3 r\sqrt{1-\frac{(1-t^2)(1-s^2)}{1-t^2(1-s^2)}}\bigg)\times\hat{\sigma}_{d-3}(|\xi|b_1 ts)\\
&\times\hat{\sigma}_{d-3}(|\eta|(\sqrt{2}/\sqrt{3})\sqrt{1-t^2(1-s^2)}r)\frac{r^{d-3}}{\sqrt{1-r^2}}\frac{s^{d-3}}{\sqrt{1-s^2}}\frac{t^{d-2}}{\sqrt{1-t^2}}\,drdsdt.\end{align}

\noindent
Let $\sin\theta :=\sin\theta(\delta,\eta) = a_3$ and $\cos\theta_1 :=\cos\theta(\xi,e_3)=b_1$. The first lemma to obtain is:\\

\noindent
{\bf Lemma 1.} For all multi-indices $\alpha,\beta,\gamma$,
\begin{equation}\label{eq E2}
|\partial^{\alpha}_{\xi}\partial^{\beta}_{\delta}\partial^{\gamma}_{\eta}\hat{\mu}(\xi,\delta,\eta)|\lesssim_{\alpha,\beta,\gamma} (1+\min(|\delta|,|\eta|)|\sin\theta|)^{-(d-3)/2}(1+|\xi||\cos\theta_1|)^{-(d-3)/2}(1+|(\xi,\delta,\eta)|)^{-(d-3)/2}.\end{equation}

\noindent
Indeed, the formula $\ref{eq A2}$ gives:
\begin{align*}\hat{\sigma}_{d-3}\bigg((\sqrt{3}/2) t|\delta|a_3 r\sqrt{1-\frac{(1-t^2)(1-s^2)}{1-t^2(1-s^2)}}\bigg) &=\frac{2\pi J_{\frac{d-4}{2}}\bigg(2\pi(\sqrt{3}/2) t|\delta||a_3| r\sqrt{1-\frac{(1-t^2)(1-s^2)}{1-t^2(1-s^2)}}\bigg)}{\bigg|(\sqrt{3}/2) t|\delta|a_3 r\sqrt{1-\frac{(1-t^2)(1-s^2)}{1-t^2(1-s^2)}}\bigg|^{\frac{d-4}{2}}}\\
\hat{\sigma}_{d-3}(|\eta|(\sqrt{2}/\sqrt{3})\sqrt{1-t^2(1-s^2)}r) &=\frac{2\pi J_{\frac{d-4}{2}}(2\pi|\eta|(\sqrt{2}/\sqrt{3})\sqrt{1-t^2(1-s^2)}r)}{||\eta|(\sqrt{2}/\sqrt{3})\sqrt{1-t^2(1-s^2)}r|^{\frac{d-4}{2}}}\\
\hat{\sigma}_{d-3}(|\xi|b_1 ts) &=\frac{2\pi J_{\frac{d-4}{2}}(2\pi|\xi||b_1| ts)}{||\xi|b_1 ts|^{\frac{d-4}{2}}}.\end{align*}
Then from {\bf Section 2}, {\it Remark 1} and in particular its $\ref{eq A1}$, one obtain that every $(\alpha,\beta,\gamma)$ partial derivatives of $\hat{\mu}$ in terms of $\xi,\delta,\eta$, respectively, can be dominated by,
\begin{multline}\label{eq E3}
C_{\alpha,\beta,\gamma}|\xi^{\nu_1}\delta^{\nu_2}\eta^{\nu_3}|(1+|\xi||b_1|)^{p_{\nu_1}}|(1+|\delta||a_3|)^{p_{\nu_2}}(1+|\eta|)^{p_{\nu_3}}\\
\times\int_0^1\int_0^1\int_0^1 \bigg|J_{u_2}\bigg(2\pi(\sqrt{3}/2) t|\delta||a_3| r\sqrt{1-\frac{(1-t^2)(1-s^2)}{1-t^2(1-s^2)}}\bigg)\bigg|\times \big|J_{u_1}(2\pi|\xi||b_1| ts)\big|\\ \times\big|J_{u_3}(2\pi|\eta|(\sqrt{2}/\sqrt{3})\sqrt{1-t^2(1-s^2)}r)\big|\\
\times\frac{1}{t^{d-4}r^{d-4}s^{(d-4)/2}\bigg(1-\frac{(1-t^2)(1-s^2)}{1-t^2(1-s^2)}\bigg)^{(d-4)/4}\bigg(1-t^2(1-s^2)\bigg)^{(d-4)/4}}\frac{r^{d-3}}{\sqrt{1-r^2}}\frac{s^{d-3}}{\sqrt{1-s^2}}\frac{t^{d-2}}{\sqrt{1-t^2}}\,drdsdt\end{multline}
where $u_{i}>0$ and $\nu_{i}$'s are multi-indices, such that,
\begin{align}\label{eq E4}
\nonumber |\xi^{\nu_1}|(1+|\xi||b_1|)^{p_{\nu_1}} &= O((1+|\xi||b_1|)^{-\frac{d-4}{2}})\\
\nonumber |\delta^{\nu_2}|(1+|\delta||a_3|)^{p_{\nu_2}} &= O((1+|\delta||a_3|))^{-\frac{d-4}{2}})\\
|\eta^{\nu_3}|(1+|\eta|)^{p_{\nu_3}} &= O((1+|\eta|)^{-\frac{d-4}{2}}).\end{align}
One also observes that in $\ref{eq E3}$:
\begin{equation*}(1-t^2(1-s^2))(1-\frac{(1-t^2)(1-s^2)}{1-t^2(1-s^2)})= s^2,\end{equation*}
thus, the non-Bessel factors in the integral in $\ref{eq E3}$ is dominated by:
\begin{equation}\label{eq E5}
C\frac{rst^2}{\sqrt{1-r^2}\sqrt{1-s^2}\sqrt{1-t^2}}\end{equation}
which is integrable on $[0,1]\times[0,1]\times[0,1]$. Lastly from {\bf Section 2}, one can deduce the following growth rate for Bessel factors in the integral $\ref{eq E3}$:
\begin{align}\label{eq E6}
\nonumber &\big|J_{u_1}(2\pi|\xi||b_1| ts)\big|\lesssim (1+|\xi||b_1|)^{-1/2}\\
\nonumber &\bigg|J_{u_2}\bigg(2\pi(\sqrt{3}/2) t|\delta||a_3| r\sqrt{1-\frac{(1-t^2)(1-s^2)}{1-t^2(1-s^2)}}\bigg)\bigg|\lesssim (1+|\delta||a_3|)^{-1/2}\\
&\big|J_{u_3}(2\pi|\eta|(\sqrt{2}/\sqrt{3})\sqrt{1-t^2(1-s^2)}r)\big|\lesssim (1+|\eta|)^{-1/2}.\end{align}
Then, due to the symmetry between $\delta,\eta$ and between $\xi,\delta,\eta$, $\ref{eq E3},\ref{eq E4},\ref{eq E5},\ref{eq E6}$ imply $\ref{eq E2}$.\\

\noindent
{\it Remark 3:} The relation between $\xi$ and $e_3$ might seem a bit artificial. However, one notes that in the calculations above, $(\xi,\delta,\eta)$ is rotated so that $\xi,\delta,\eta$ are all in the $3$-space created by $e_1,e_2,e_3$ and $\delta,\eta$ are both in the $2$-plane spanned by $e_1,e_2$. Hence the relation between $\xi$ and $e_3$, expressed by $\cos\theta_1$, can be thought of as the relation between $\xi$ and $span(\delta,\eta)$ or the relation between $\xi$ and $\delta,\eta$.\\

\noindent
{\it Remark 4:} One can observe from $\ref{eq S24}$ and $\ref{eq E2}$ that when any of the vectors $\xi,\delta,\eta$ is the zero vector or when $a_3=0$ or $b_1=0$ the derived decay for the derivatives has a worse bound than what is shown in $\ref{eq E2}$. It's interesting since it might hint as a sharp change in behavior for the derivatives. It's also non-interesting for the author as these rare cases happen with zero probability.

\section{Decomposition of $\hat{\mu}$}

\subsection{Partition in $|(\eta,\delta)|$}

\noindent
Let $\Phi_0\in C_{c}^{\infty}(\mathbb{R}^{2d})$ satisfying $\chi_{\mathbb{B}^{2d}}\leq\Phi_0\leq\chi_{2\mathbb{B}^{2d}}$. For $i\geq 1$, define \begin{equation*}\Phi_{i}((\eta,\delta))=\Phi_0(2^{-i}(\eta,\delta))-\Phi_0(2^{-i+1}(\eta,\delta)).\end{equation*} Then define 
\begin{equation*}\phi_{i}(\xi,\delta,\eta)=1(\xi)\Phi_{i}(\eta,\delta).\end{equation*}
Note that $\phi_{i}$ is supported on the strip $\{(\xi,\delta,\eta): 2^{i-1}\leq|(\eta,\delta)|\leq 2^{i+1}\}$ and $\sum_{i\geq 0}\phi_{i}(\xi,\delta,\eta)=1.$ 

\subsection{Partition in $|\eta|/|\delta|$}

\noindent
Let $\epsilon$ be a sufficiently small positive quantity. Let $\Psi\in C_{c}^{\infty}(\mathbb{R})$ such that $\chi_{[0,1]}\leq\Psi\leq\chi_{[-\epsilon,1+\epsilon]}$. For $i\in\mathbb{Z}$, define
\begin{equation*}\Psi_{i}(t) = \frac{\Psi(t-i)}{\sum_{k}\Psi(t-k)}.\end{equation*}
Note that $\Psi_{i}$ is defined on the interval $[-\epsilon+j, j+1+\epsilon]$. For $i\geq 0$, define
\begin{equation*}\psi_{i}((\xi,\delta,\eta))=1(\xi)(\Psi_{i}(\log|\eta|-\log|\delta|)+\Psi_{-i-1}(\log|\eta|-\log|\delta|)),\end{equation*}
which is supported on $\{(\xi,\delta,\eta):2^{-\epsilon}2^{-i}\leq\frac{\min(|\eta|,|\delta|)}{\max(|\eta|,|\delta|)}\leq 2^{\epsilon}2^{-i+1}\}$. Also for $i\geq 0$, define
\begin{equation*}\psi^{i}((\xi,\delta,\eta))=\sum_{k\geq i}\psi_{i}((\xi,\delta,\eta)).\end{equation*}
Then $\psi^{i}$ is supported on $\{(\xi,\delta,\eta):\frac{\min(|\eta|,|\delta|)}{\max(|\eta|,|\delta|)}\leq 2^{\epsilon}2^{-i+1}\}$ and is identically $1$ on $\{(\xi,\delta,\eta):\frac{\min(|\eta|,|\delta|)}{\max(|\eta|,|\delta|)}\leq 2^{-\epsilon}2^{-i}\}$.\\

\noindent
{\it Remark 5:} To see the second claim above, let $t=\log|\eta|-\log|\delta|$ and consider \begin{equation*}\sum_{j\geq i}\Psi_{j}(t) + \Psi_{-j-1}(t) =\frac{\sum_{j\geq i}\Psi((\log|\eta|-\log|\delta|)-j)}{\sum_{k}\Psi((\log|\eta|-\log|\delta|)-k)} + \frac{\sum_{j\geq i}\Psi((\log|\eta|-\log|\delta|)+j+1)}{\sum_{k}\Psi((\log|\eta|-\log|\delta|)-k)}.\end{equation*}
Suppose $2^{-\epsilon}2^{-i}|\eta|\geq|\delta|$. Then the second sum on the RHS above is not "activated" (being zero), and similarly, the bottom sum is reduced to $\sum_{k\geq i}\Psi((\log|\eta|-\log|\delta)-k)$. One is left with $(\sum_{j\geq i}\Psi((\log|\eta|-\log|\delta|)-j))/(\sum_{j\geq i}\Psi((\log|\eta|-\log|\delta|)-j))=1$.

\subsection{Partition in $|\sin\theta|$}

\noindent
Let $\mathcal{P}\in C_{c}^{\infty}(\mathbb{R})$ be such that $\chi_{[-1,1]}\leq\mathcal{P}\leq\chi_{[-2,2]}$. For $i\in\mathbb{Z}$, define 
\begin{equation*}\mathcal{P}_{i}(t) = \mathcal{P}(2^{-2i}t)-\mathcal{P}(2^{-2(i-1)}t).\end{equation*} 
Then for $i\geq 1$, define
\begin{equation*}\rho_{i}((\xi,\delta,\eta))=1(\xi)\mathcal{P}_{-i}\bigg(1-\bigg(\frac{\eta\cdot\delta}{|\eta||\delta|}\bigg)^2\bigg),\end{equation*}
then $\rho_{i}$ is supported on $\{(\xi,\delta,\eta): 2^{-(i+1)}\leq|\sin\theta|\leq 2^{-(i-1)}\}$. Also, define
\begin{equation*}\rho^{i}((\xi,\delta,\eta))=1(\xi)\left(\sum_{j\geq i}\mathcal{P}_{-j}\bigg(1-\bigg(\frac{\eta\cdot\delta}{|\eta||\delta|}\bigg)^2\bigg)\right).\end{equation*}
This means that $\rho^{i}$ is supported on $\{(\xi,\delta,\eta):|\sin\theta|\leq 2^{-(i-1)}\}$ and is identically $1$ on $\{(\xi,\delta,\eta):|\sin\theta|\leq 2^{-(i+1)}\}$. For $i=0$, define
\begin{equation*}\rho_0((\xi,\delta,\eta))=1(\xi)\left(\sum_{j\leq 0}\mathcal{P}_{-j}\bigg(1-\bigg(\frac{\eta\cdot\delta}{|\eta||\delta|}\right)^2\bigg)\bigg).\end{equation*}
Then $\rho_0((\xi,\delta,\eta))$ is supported on $\{(\xi,\delta,\eta): |\sin\theta|\geq 1/2\}$.

\subsection{Partition in $|\xi|$}

\noindent
Let $\mathcal{Z}_0\in C^{\infty}_{c}(\mathbb{R}^{d})$ satisfying $\chi_{\mathbb{B}^{d}}\leq\mathcal{Z}_0\leq\chi_{2\mathbb{B}^{d}}$. For $i\geq 1$, define \begin{equation*}\mathcal{Z}_{i}(\xi) = \mathcal{Z}_0(2^{-i}(\xi))-\mathcal{Z}_0(2^{-i+1}(\xi)).\end{equation*}
Then for $i\geq 0$, define
\begin{equation*}\zeta_{i}((\xi,\delta,\eta))=\mathcal{Z}_{i}(\xi)1((\eta,\delta)).\end{equation*}
Certainly, $\zeta_{i}$ is supported on $\{(\xi,\delta,\eta): 2^{i-1}\leq|\xi|\leq 2^{i+1}\}$ and $\sum_{i\geq 0}\zeta_{i}((\xi,\delta,\eta))=1$. 

\subsection{Partition in $|\cos\theta_1|$}

\noindent
Let $\mathcal{P}\in C_{c}^{\infty}(\mathbb{R})$ and $\mathcal{P}_{i}$ be as before. Then for $i\geq 1$, define
\begin{equation*}\rho_{1,i}((\xi,\delta,\eta))=1(\xi)\mathcal{P}_{-i}\bigg(\bigg(\frac{\xi\cdot e_3}{|\xi|}\bigg)^2\bigg),\end{equation*}
then $\rho_{1,i}$ is supported on $\{(\xi,\delta,\eta): 2^{-(i+1)}\leq|\cos\theta_1|\leq 2^{-(i-1)}\}$. Also, define
\begin{equation*}\rho_1^{i}((\xi,\delta,\eta))=1(\xi)\left(\sum_{j\geq i}\mathcal{P}_{-j}\bigg(\bigg(\frac{\xi\cdot e_3}{|\xi|}\right)^2\bigg)\bigg).\end{equation*}
This means that $\rho_1^{i}$ is supported on $\{(\xi,\delta,\eta):|\cos\theta_1|\leq 2^{-(i-1)}\}$ and is identically $1$ on $\{(\xi,\delta,\eta):|\cos\theta_1|\leq 2^{-(i+1)}\}$. For $i=0$, define
\begin{equation*}\rho_{1,0}((\xi,\delta,\eta))=1(\xi)\left(\sum_{j\leq 0}\mathcal{P}_{-j}\bigg(\bigg(\frac{\xi\cdot e_3}{|\xi|}\right)^2\bigg)\bigg).\end{equation*}
Then $\rho_{1,0}((\xi,\delta,\eta))$ is supported on $\{(\xi,\delta,\eta): |\cos\theta_1|\geq 1/2\}$.

\subsection{The decomposition}

\noindent
Let $m=\hat{\mu}$. Fix $n\in\mathbb{N}_0$. If $i=0$, define $m_{0,0,0,0} = m\phi_0\psi^0\zeta_0$.\\

\noindent
When $1\leq i$, $0\leq j\leq i-1$, $0\leq k<\lfloor\frac{i-j}{2}\rfloor$, define,
\begin{equation*}m_{i,j,k}=
m\phi_{i}\psi_{j}\rho_{k}\rho_{1,k}\zeta_{i}.\end{equation*}
When $0\leq j\leq i-1$, $k=\lfloor\frac{i-j}{2}\rfloor$, define,
\begin{equation*}m_{i,j,\lfloor\frac{i-j}{2}\rfloor}=
m\phi_{i}\psi_{j}\rho^{\lfloor\frac{i-j}{2}\rfloor}\rho_1^{\lfloor\frac{i-j}{2}\rfloor}\zeta_{i}.\end{equation*}
When $i=j\geq 1$, define
\begin{equation*}m_{i,i,0}=
m\phi_{i}\psi^{i}\zeta_{i}.\end{equation*}

\noindent
{\it Remark 6:} The decomposition above is inspired by what was done in [7]. The philosophy behind this decomposition scheme is as follows. Starting with a "level" $R=|(\xi,\delta,\eta)|$, one wants $|\xi|\sim |(\delta,\eta)|$ and the ratio $|\eta|/|\delta|$ (or $|\delta|/|\eta|$) as well as the values of $|\sin\theta|,|\cos\theta_1|$ to be divided into possible ranges. When the level $R$ is small, there is no need to divide up cases for the angles. Care needs to be taken the "border" cases - either $j=i$ or $k=\lfloor\frac{i-j}{2}\rfloor$. For the "non-typical" cases, the saving grace is either these cases are finite in number or the derivatives of the considered functions are vanishing in the appropriate regions - by construction, derivatives of $\psi^{i},\rho^{i},\rho_1^{i}$ are vanishing on $\{(\xi,\delta,\eta):\frac{\min(|\eta|,|\delta|)}{\max(|\eta|,|\delta|)}\leq 2^{-\epsilon}2^{-i}\}$, $\{(\xi,\delta,\eta):|\sin\theta|\leq 2^{-(i+1)}\}$ and $\{(\xi,\delta,\eta):|\cos\theta_1|\leq 2^{-(i+1)}\}$, respectively.\\

\noindent
For the "typical" cases, which are when $i\geq j$, $i-j\geq 2k$, in the support of $m_{i,j,k}$ one has the following:
\begin{align}\label{eq D1} 
\nonumber 2^{2(i-1)} &\leq|\eta|^2+|\delta|^2\leq 2^{2(i+1)}\\
\nonumber 2^{-\epsilon}2^{-j} &\leq\frac{\min(|\eta|,|\delta|)}{\max(|\eta|,|\delta|)}\leq 2^{-\epsilon}2^{-j+1}\\
\nonumber 2^{-(k+1)} &\leq|\sin\theta|\leq 2^{-(k-1)}\\
\nonumber 2^{-(k+1)} &\leq|\cos\theta_1|\leq 2^{-(k-1)}\\
2^{i-1} &\leq |\xi|\leq 2^{i+1},\end{align}
which implies $2^{i+2}\geq|\eta|+|\delta|\geq 2^{i-1}$, and hence
\begin{align}\label{eq D2}
\nonumber 2^{i-j+2}\geq\min(|\eta|,|\delta|) &\geq 2^{i-j-3}\\
|(\xi,\delta,\eta)|\geq (1/2)(|\xi|+|(\delta,\eta)|) &\geq 2^{i-1}.\end{align}
Then, from these calculations and the definitions of $\phi_{i},\psi_{j},\rho_{k},\rho_{1,k}\zeta_{i}$, one has that in the support of $m_{i,j,k}$, 
\begin{align}\label{eq D3}
\nonumber\|\partial^{\alpha}_{\xi}\partial^{\beta}_{\delta}\partial^{\gamma}_{\eta}\phi_{i}\|_{L^{\infty}},\|\partial^{\alpha}_{\xi}\partial^{\beta}_{\eta}\partial^{\gamma}_{\delta}\zeta_{i}\|_{L^{\infty}} &\leq C\\
\nonumber\|\partial^{\alpha}_{\xi}\partial^{\beta}_{\delta}\partial^{\gamma}_{\eta}\psi_{j}\|_{L^{\infty}} &\lesssim 2^{-(|\beta|+|\gamma|)(i-j)}\\
\nonumber \|\partial^{\alpha}_{\xi}\partial^{\beta}_{\delta}\partial^{\gamma}_{\eta}\rho_{k}\|_{L^{\infty}} &\lesssim 2^{-(|\beta|+|\gamma|)(i-j-2k)}\\
\|\partial^{\alpha}_{\xi}\partial^{\beta}_{\delta}\partial^{\gamma}_{\eta}\rho_{1,k}\|_{L^{\infty}} &\lesssim 2^{-(|\alpha|)(n-2k)}.\end{align}
Finally, from definitions and from $\ref{eq D1},\ref{eq D2}$, the support of $m_{i,j,k,n}$, for all cases except for when $k=0$ or $i=j$, is contained in,
\begin{equation*}\{(\xi,\delta,\eta): \min(|\eta|,|\delta|)\leq 2^{i-j+2}; |\xi|,\max(|\eta|,|\delta|)\leq 2^{i+1};|\sin\theta|,|\cos\theta_1|\leq 2^{-(k-1)}\}.\end{equation*}

\noindent
Temporarily assume $|\xi|=|\delta|=|\eta|=1$ for a moment. Then the volume of the region $\{(\xi,\delta,\eta): |\sin\theta|,|\cos\theta_1|\leq 2^{-(k-1)}\}$ is,
\begin{align*}&\left(\int_{\mathbb{S}^{d-1}}\int_{\mathbb{S}^{d-1}} \chi_{\{|\sin\theta(\omega_1,\omega_2)|\leq 2^{-k+1}\}}(\omega_1,\omega_2)\,d\sigma(\omega_1)d\sigma(\omega_2)\right)\times \left(\int_{\mathbb{S}^{d-1}} \chi_{\{|\cos\theta(\omega_3,e_3)|\leq 2^{-k+1}\}}(\omega_3)\,d\sigma(\omega_3)\right)\\
&=_{d}\left(\int_{\mathbb{S}^{d-1}} \chi_{\{|\sin\theta(\omega_1,e_1)|\leq 2^{-k+1}\}}(\omega_1)\,d\sigma(\omega_1)\right)
\times \left(\int_{\mathbb{S}^{d-1}} \chi_{\{|\cos\theta(\omega_3,e_3)|\leq 2^{-k+1}\}}(\omega_3)\,d\sigma(\omega_3)\right)\\
&=_{d}\left(\int_{\mathbb{S}^{d-1}} \chi_{\{|\cos\theta(\omega_1,e_1)|\leq 2^{-k+1}\}}(\omega_1)\,d\sigma(\omega_1)\right)^2\end{align*}
by symmetry and by the fact that $\sin(\cdot)$ and $\cos(\cdot)$ are $\pi/2$ off-phase versions of each other.

\noindent
That means such support volume is bounded by,
\begin{align}\label{eq D4}
\nonumber &\int_0^{2^{i+1}}\int_0^{2^{i+1}}\int_0^{2^{i-j+2}}\left(\int_{\mathbb{S}^{d-1}} \chi_{\{|\cos\theta(\omega_1,e_1)|\leq 2^{-k+1}\}}(\omega_1)\,d\sigma(\omega_1)\right)^2\,drdsdt\\
\nonumber =_{d} & 2^{2di}2^{d(i-j)}\left(\int_{\mathbb{S}^{d-1}} \chi_{\{|\cos\theta(\omega_1,e_1)|\leq 2^{-k+1}\}}(\omega_1)\,d\sigma(\omega_1)\right)^2\\
\nonumber =_{d} & 2^{2di}2^{d(i-j)}\left(\int_{\mathbb{S}^{d-1}}\chi_{\{|\omega\cdot e_1|\leq 2^{-k+1}\}}\,d\sigma(\omega_1)\right)^2\\
\nonumber =_{d} & 2^{2di}2^{d(i-j)}\left(\int_0^{2^{-k+1}}\frac{r^{d-2}}{\sqrt{1-r^2}}\,dr\right)^2\\
\lesssim_{d} &2^{2di}2^{d(i-j)}2^{-2k(d-3)},\end{align}
where one uses the slicing formula in the second to last step.\\ 

\noindent
{\it Remark 7:} In the cases when $k=0$ or $i=j$, then by definitions, there is no need for the two spherical integrations in $\ref{eq D4}$. One is left with the dominants $O_{d}(2^{2di}2^{d(i-j)})$ or $O_{d}(2^{2di})$ instead. 

\section{$L^2$ boundedness}

\noindent
Let 
\begin{equation*} m_{i}=\sum_{j=0}^{i}\sum_{k=0}^{\lfloor\frac{i-j}{2}\rfloor}m_{i,j,k}.\end{equation*} 
Then let $T_{i}, T_{i,j,k}$ denote the operator associated with the multipliers $m_{i}, m_{i,j,k}$, respectively. In other words, define
\begin{align*}T_{i}(f,g,h)(x) &= \mathcal{F}^{-1}(m_{i}\hat{f}\otimes\hat{g}\otimes\hat{h})(x,x,x)\\
T_{i,j,k,n}(f,g,h)(x) &= \mathcal{F}^{-1}(m_{i,j,k}\hat{f}\otimes\hat{g}\otimes\hat{h})(x,x,x)\end{align*}
with $f,g,h\in\mathcal{S}(\mathbb{R}^{d})$. One has that, 
\begin{equation}\label{eq B1}
\|T\|_{L^2\times L^2\times L^2\to L^{2/3}}\leq\sum_{i\geq 0}\|T_{i}\|_{L^2\times L^2\times L^2\to L^{2/3}}\leq\sum_{i\geq 0}\left(\sum_{j=0}^{i}\sum_{k=0}^{\lfloor\frac{i-j}{2}\rfloor}\|T_{i,j,k}\|_{L^2\times L^2\times L^2\to L^{2/3}}\right).\end{equation}
Hence one hopes that the sum on the farthest left is summable in $i$.\\

\noindent
Now $\ref{eq D3}$ and the choices of $j,k$ implies that within its support, the derivatives of $m_{i,j,k}$ are uniformly bounded in $i,j,k$. That means, from definitions of $\phi_{i},\psi_{j},\rho_{k},\rho_{1,k},\zeta_{i}$, $\ref{eq E2},\ref{eq D1},\ref{eq D2}$, one has,
\begin{equation}\label{eq B2}
\|\partial^{\alpha}_{\xi}\partial^{\beta}_{\eta}\partial^{\gamma}_{\delta}m_{i,j,k}(\xi,\eta,\delta)\|_{L^{\infty}}\lesssim_{d,\alpha,\beta,\gamma}2^{-2i(d-3)}2^{j(d-3)/2}2^{k(d-3)}.\end{equation}
One simply puts all the findings in $\ref{eq D4}, \ref{eq B2}$ and {\it Remark 7} together, and invoke {\bf Lemma A} and its {\it Remark 2}, to arrive at,
\begin{equation}\label{eq B3}
\|T_{i,j,k,n}\|_{L^2\times L^2\times L^2\to L^{2/3}}\lesssim_{d}2^{-(2i/3)(d-3)}2^{j(d-3)/2}2^{k(d-3)} \cdot(2^{2di}2^{d(i-j)}2^{-2k(d-3)})^{1/3}.
\end{equation}
Summing $\ref{eq B3}$ in terms of $j,k$ gives the following dominant
\begin{align*} 
\sum_{j=0}^{i}\sum_{k=0}^{\lfloor\frac{i-j}{2}\rfloor}\|T_{i,j,k}\|_{L^2\times L^2\times L^2\to L^{2/3}}&\lesssim_{d}2^{i(9/2-d/2)}\sum_{j=0}^{i}2^{jd/6}2^{-3j/2}\sum_{k=0}^{\lfloor\frac{i-j}{2}\rfloor}2^{k(d-3)/3}\\
&\lesssim_{d} 2^{i(9/2-d/2)}2^{id/6}2^{-3i/2}2^{i(d-3)/6}\\
&\lesssim_{d} 2^{-id/6+i5/2}.\end{align*}
Hence, in order for the last sum in $\ref{eq B1}$ to be summable in $i$, one must have, $d> 15$. This establishes the dimensional sufficiency for the boundedness of $T: L^2\times L^2\times L^2\to L^{2/3}$.

\section{Bounded region}

\noindent
Recall that $T(f,g,h)$ can be written in the following form:
\begin{multline*}
T(f,g,h)(x) =\int_{SO(d)}f(x-Re_1) g(x- [R(1/2)e_1+R(\sqrt{3}/2)e_2])\\
\times h(x-[R(1/2)e_1+R(1/2\sqrt{3})e_2+R(\sqrt{2}/\sqrt{3})e_3])\,dR.\end{multline*}
If this times, one precomposes $R\in SO(d)$ with a rotation whose transpose fixes the two vectors $Re_1$, $R(1/2)e_1+R(\sqrt{3}/2)e_2$, then one obtains the following form of $T(f,g,h)$:
\begin{multline}\label{eq R1}
T(f,g,h)(x) =\int_{SO(d)}f(x-Re_1) g(x- R[(1/2)e_1+(\sqrt{3}/2)e_2])\\
\times\int_{SO(d-2)}h(x-[RR'(1/2)e_1+RR'(1/2\sqrt{3})e_2+RR'(\sqrt{2}/\sqrt{3})e_3])\,dR'dR.\end{multline}
An application of the Minkowski's integral inequality [8] to $\ref{eq R1}$ gives,
\begin{multline}\label{eq R2}
|T(f,g,h)(x)|\leq\int_{SO(d-2)}\int_{SO(d)}|f(x-Re_1)|\big|g(x- [R(1/2)e_1+R(\sqrt{3}/2)e_2])\big|\\
\times \big|h(x-[RR'(1/2)e_1+RR'(1/2\sqrt{3})e_2+RR'(\sqrt{2}/\sqrt{3})e_3])\big|\,dRdR'.\end{multline}
If one recognizes that,
\begin{equation*}\Delta(f,g)(x)=\int_{SO(d)}f(x-Re_1)g(x- [R(1/2)e_1+R(\sqrt{3}/2)e_2])\,dR\end{equation*}
is the triangle operator considered in [7], one obtains from $\ref{eq R2}$
\begin{equation}\label{eq R3}
|T(f,g,h)(x)|\lesssim_{d}\|h\|_{\infty}|\Delta(|f|,|g|)(x)|.\end{equation}
Obtained in [7] is an improved range for boundedness of the operator $\Delta$, when $d>5$, which the found requirement $d>15$ satisfies. That range $\Delta: L^{p}\times L^{q}\to L^{r}$ is as follows
\begin{align}\label{eq R4}
\nonumber (1/p,1/q) &\in Conv((0,1),(1,0),(0,0),((3d-2)/5d,(3d-2)/5))\\
1/r &= 1/p+1/q.
\end{align}
Hence, using multilinear interpolation theory [1], $\ref{eq R3},\ref{eq R4}$ and the fact that $T(f,g,h)$ is symmetric in $f,g,h$ implies that $T: L^{p}\times L^{q}\times L^{s}\to L^{r}$ is bounded in:
\begin{align}\label{eq R5}
\nonumber (1/p,1/q,1/s) &\in Conv(S)\\
1/r &= 1/p+1/q+1/s\end{align}
where $S=\{(1,0,0),(0,1,0),(0,0,1),(0,0,0),((3d-2)/5d,(3d-2)/5d,0),(0,(3d-2)/5d,(3d-2)/5d),((3d-2)/5d,0,(3d-2)/5d))\}$. Note that this range includes the mentioned Banach range. However this convex hull $Conv(S)$ does not include the point $(1/2,1/2,1/2)$ - nor does the Banach range. As shown in the previous sections $T: L^2\times L^2\times L^2\to L^{2/3}$ as long as $d>15$. To see that indeed $(1/2,1/2,1/2)\not\in Conv(S)$, let $p_0=\frac{5d}{3d-2}$ and suppose that there are $0<t_{i}<1$ such that $\sum_{i=1}^6 t_{i}=1$ and
\begin{align}\label{eq R6}
\nonumber t_1+t_2+t_3+ &t_4+t_5+t_6=1\\
t_1\frac{1}{p_0}(e_1+e_2)+t_1\frac{1}{p_0}(e_1+e_2)+t_1\frac{1}{p_0}(e_1+e_2)+ &t_4e_1+t_5e_2+t_6e_3=(1/2,1/2,1/2).\end{align}
Collecting like terms in $\ref{eq R6}$ yields the following system:
\begin{align}\label{eq R7}
\nonumber t_1+t_2+t_3+t_4+t_5+t_6 &=1\\
\nonumber t_4+t_1/p_0+t_3/p_0 &=1/2\\
\nonumber t_5+t_2/p_0+t_1/p_0 &=1/2\\
t_6+t_2/p_0+t_3/p_0 &=1/2.\end{align}
Adding the last three equations in $\ref{eq R7}$ yields another system
\begin{align*}
t_1+t_2+t_3+t_4+t_5+t_6 &=1\\
(t_4+t_5+t_6)+(2/p_0)(t_1+t_2+t_3) &=3/2.\end{align*}
Thus it's enough to find if there are two numbers $0<t,s<1$ such that $t+(2/p_0)s=3/2$ and $t+s=1$. But this gives a linear system whose solution $(t,s)$ implies that $t=\frac{5d}{2d-8}>1$: a contradiction.\\

\noindent
But this means that one can further improves the range in $\ref{eq R5}$. Hence utilizing the multilinear interpolation theory again, one concludes that $T: L^{p}\times L^{q}\times L^{s}\to L^{r}$ in:
\begin{align}\label{eq R8}
\nonumber (1/p,1/q,1/s) &\in Conv(S')\\
1/r &= 1/p+1/q+1/s\end{align}
where $S'=\{(1,0,0),(0,1,0),(0,0,1),(0,0,0),((3d-2)/5d,(3d-2)/5d,0),(0,(3d-2)/5d,(3d-2)/5d),((3d-2)/5d,0,(3d-2)/5d),(1/2,1/2,1/2)\}$.


\end{document}